\magnification=\magstephalf
{\centerline{ \bf THE \lq\lq NORTH POLE PROBLEM\rq\rq \/ AND RANDOM ORTHOGONAL MATRICES}}
\bigskip\bigskip
\noindent{{MORRIS L. EATON}

\noindent{School of Statistics, University of Minnesota,  224 Church St SE, Minneapolis, 
MN 55455} 

\noindent{email: eaton@stat.umn.edu}

\bigskip
\noindent {ROBB J. MUIRHEAD}

\noindent {Statistical Research and 
Consulting Center, Pfizer Inc.,  50 Pequot Ave, New 
London, CT 06320}

\noindent{email: robb.j.muirhead@pfizer.com}



\bigskip\bigskip

\noindent{\it Abstract}
\medskip
\noindent This paper is motivated by the following observation. Take a $3\times 3$ random (Haar distributed)  orthogonal matrix $\Gamma$, and use it to \lq\lq rotate" the north pole, $x_0$ say, on the unit sphere in $R^3$. This then gives a point $u=\Gamma x_0$ that is uniformly distributed on the unit sphere. Now use the same orthogonal matrix to transform $u$, giving $v=\Gamma u=\Gamma^2 x_0$. Simulations reported in Marzetta et al (2002) suggest that $v$ is more likely to be in the northern hemisphere than in the southern hemisphere,   and, morever, that $w=\Gamma^3 x_0$ has higher probability of being closer to the poles $\pm x_0$ than the uniformly distributed point $u$.  In this paper we prove these results, in the general setting of dimension $p\ge 3$, by deriving the exact distributions of the relevant components of $u$ and $v$.  The essential questions answered are the following. Let $x$ be any fixed point on the unit sphere in $R^p$, where $p\ge 3$. What are the distributions of $U_2=x'\Gamma^2 x$ and $U_3=x'\Gamma^3 x$? It is clear by orthogonal invariance that these distribution do not depend on $x$, so that we can, without loss of generality, take $x$ to be $x_0=(1,0,...,0)'\in R^p$. Call this the \lq\lq north pole\rq\rq. Then $x_0'\Gamma^ k x_0$ is the first component of the vector $\Gamma^k x_0$. We derive stochastic representations for the exact distributions of $U_2$ and $U_3$ in terms of random variables with known distributions.

\bigskip
\noindent {\bf 1. Introduction}
\medskip 
\noindent Random orthogonal matrices have long played an important role in multivariate statistical analysis (see, for example, James (1954), Muirhead (1982)); more recently such matrices have received attention in the encryption and signal processing literature (see Sloane (1983), Marzetta et al (2002)). Although the uniform (Haar) distribution on the group ${\cal O}_p$  of $p\times p$  orthogonal matrices is fairly well understood, questions regarding powers of uniformly distributed random orthogonal matrices have arisen in the past few years. For ${\cal O}_p$, as well as for other classical compact Lie groups, Rains (1997) provides an intriguing description of the distribution of high powers of Haar distributed group elements. Some of the results in Rains (1997) are reproduced in Marzetta et al (2002) where it is demonstrated via numerically constructed graphic displays what happens on ${\cal O}_3$ for low powers of Haar distributed orthogonals. For example, suppose $\Gamma$ is Haar distributed on ${\cal O}_3$ and $x_0=(1,0,0)'$ is a unit vector in $R^3$. Figure 3 of Marzetta et al (2002, p. 947)) shows the curious result that simulated observations on $\Gamma^2 x_0$  bunch up around $x_0$  (compared to the uniformly distributed vector $\Gamma x_0$), while simulated observations on $\Gamma^3 x_0$ bunch up around $\pm x_0$. This phenomenon is an example of what we call the \lq\lq north pole problem\rq\rq.
\smallskip
   To justify our \lq\lq north pole\rq\rq nomenclature and describe the results below, suppose that $\Gamma$ is uniformly distributed on ${\cal O}_p$, where $p\ge 3$. Consider a column vector $x\in R^p$ with $x'x=\|x\|^2=1$.  For a positive integer $k$, let
$$U_k=x'\Gamma^k x.\eqno(1.1)$$
Obviously $U_k$ has a distribution that depends on $p$, but we will suppress this dependence in our notation. 
 In what follows, we will use the \lq\lq law\rq\rq  notation, ${\cal L}(\cdot)$, to denote the distribution or probability law of \lq\lq\ $\cdot$\ \rq\rq, where \lq\lq\ $\cdot$\ \rq\rq can be a random variable, a random vector, or a random matrix. In this notation, the uniform distribution on ${\cal O}_p$ is characterized by its invariance, namely
$${\cal L}(\Gamma)={\cal L}(g\Gamma)= {\cal L}(\Gamma h),\eqno(1.2)$$
 for all for all $g,h\in{\cal O}_p$.  Equation (1.2) immediately implies 
$${\cal L}(x'\Gamma^k x)= {\cal L}(y'\Gamma^k y)\eqno(1.3)$$ for all $x,y\in R^p$ with $\|x\|=\|y\|=1$. Therefore the distribution of $U_k$ in (1.1) does not depend on $x$, so in what follows we simply take $$x=x_0=\pmatrix{1\cr
0\cr
\vdots\cr
0\cr}.\eqno(1.4)$$
 We will then arbitrarily call $x_0$ the \lq\lq north pole\rq\rq on the unit sphere ${\cal S}_p=\{x\in R^p\vert\|x\|=1\}$. Next, observe that $x_0x_0'(\Gamma^k x_0)$ is the orthogonal projection of $\Gamma^k x_0$ onto the 1-dimensional subspace span$\{x_0\}$.  Therefore the random variable $\vert U_k\vert$ is just the length of the orthogonal projection of $\Gamma^k x_0$ onto the span of the north pole. For this reason, we refer to the problem of finding ${\cal L}(U_k)$ as the \lq\lq north pole problem\rq\rq.
\smallskip

   The results in this paper give the exact distribution of $U_k$  for $k = 2, 3$ and for all $p\ge 3$.  Our results are expressed in terms of random variables with known distributions, and this facilitates the simulation of ${\cal L} (U_k)$  for $p$ fixed. For example, when $k = 1$, $U_1$ is the (1, 1) element of $\Gamma$ and this is well known to have the probability density on $R^1$ given by
$$f(x|p)={\Gamma({1\over 2}p)\over \Gamma({1\over 2})\Gamma({1\over 2}(p-1))}(1-x^2)^{(p-3)/2},\qquad \vert x\vert <1.\eqno(1.5)$$
Of course, $U_1^2$ has the beta distribution ${\rm Beta}({1\over 2},{1\over 2}(p-1))$ with density function
$$g_0(y)= {\Gamma({1\over 2}p)\over \Gamma({1\over 2})\Gamma({1\over 2}(p-1))}y^{{1\over 2}-1}(1-y)^{{1\over 2}(p-1)-1},\qquad 0<y<1.\eqno(1.6)$$
\smallskip
   Our approach to finding the distributions of $U_2$ and $U_3$ is based on some recent results in Eaton and Muirhead (2008). These results are expressed in terms of a decomposition of the Haar distribution corresponding to the partition 
$$\Gamma=\pmatrix{\Gamma_{11}&\Gamma_{12}\cr
\Gamma_{21}&\Gamma_{22}\cr}\eqno(1.7)$$
 where $\Gamma_{11}$ is $1\times 1$,  $\Gamma_{12}$ is $1\times (p-1)$, $\Gamma_{21}$ is $(p-1)\times 1$ and $\Gamma_{22}$ is $(p-1)\times (p-1)$.
In particular, ${\cal L}(\Gamma_{11})$ is specified by the density (1.5). Eaton and Muirhead (2008) give the conditional distribution of $(\Gamma_{21},\Gamma_{12})$ given $\Gamma_{11}$ and the conditional distribution of  $\Gamma_{22}$ given $(\Gamma_{11},\Gamma_{21},\Gamma_{12})$. These results are summarized in Section 2 of this paper.
\smallskip
   The results in Section 2 are then used in Section 3 to give random variable presentations of ${\cal L}(U_2)$ and ${\cal L}(U_3)$.  For example, here is the result for $U_2$.
 \medskip  
\noindent{\bf Theorem 1.1:}\ \  Let $\xi_1$  and $\xi_2$ be independent random variables, where $\xi_1$ has density $f(\cdot|p)$ given in (1.5) and $\xi_2$ has density $f(\cdot|p-1)$. Then
$${\cal L}(U_2)={\cal L}\bigl(\xi_1^2+(1-\xi_1^2)\xi_2\bigr).\eqno(1.8)$$

\smallskip

   The representation (1.8), which is derived in Section 3,  shows that $U_2$  can be thought of as a random convex combination (with weights $\xi_1^2$ and $(1-\xi_1^2)$) of a random variable that is 1 and the random variable $\xi_2$. This helps explain the Marzetta et al (2002) simulation result that, when $p=3$,  $U_2$  exhibits some bunching near 1.   Here are values of $P(U_2>0)$ for various dimensions $p$:

\bigskip
\settabs 11\columns
\+&Dimension $p$&&3&4&5&10&20&50&100&500\cr
\smallskip
\+&$P(U_2>0)$ &&0.71&0.68&0.66&0.62&0.59&0.56&0.54&0.52\cr

\bigskip
\noindent These probabilities suggest convergence to $1/2$ as $p\to\infty$. In fact, it is straightforward to show that ${\cal L}({\sqrt{p}}\/U_2)\to N(0,1)$ as $p\to\infty$.

\smallskip 
The distribution $U_3$ is somewhat more complicated than that of $U_2$, but is still expressible in terms of independent random variables. It is given in Section 3. 
   
\smallskip

   Powers of $\Gamma$ have received attention in the literature, especially concerning the distribution of the eigenvalues of $\Gamma^k$, $k=1,2,...$. For example, see Diaconis and Shahshahani (1994), Rains (1997), Marzetta et al (2002), and Diaconis (2003). 
 
\bigskip
  
\noindent {\bf 2.  A Decomposition Theorem}
\medskip
\noindent Let $\Gamma$ have the Haar distribution on ${\cal O}_p$ and partition  $\Gamma$ as in (1.7). In this section, we describe a recent result of Eaton and Muirhead (2008) that characterizes ${\cal L}(\Gamma)$. The basic idea is to specify\smallskip
$$\eqalign{(i)\quad &{\cal L}(\Gamma_{11})\cr
(ii)\quad &{\cal L}(\Gamma_{21},\Gamma_{12}|\Gamma_{11})\cr
(iii)\quad &{\cal L}(\Gamma_{22}|\Gamma_{11},\Gamma_{21},\Gamma_{12})\cr}\eqno(2.1)$$
\smallskip
\noindent where $(ii)$ and $(iii)$ in (2.1) refer to obvious conditional distributions. It is clear that ${\cal L}(\Gamma)$ is completely determined by $(i)$, $(ii)$ and $(iii)$ in (2.1). 
\smallskip
   In order to present a rigorous statement concerning the above conditional distributions, it is useful to modify the sample space for $\Gamma$ slightly. Let
$${\cal O}_p^+=\{g|g\in{\cal O}_p,g_{11}\in (-1,1)\},\eqno(2.2)$$
 where $g_{11}$ is the (1,1) element of the matrix $g$. Note that ${\cal O}_p^+$  has probability one under the Haar distribution on ${\cal O}_p$.
\smallskip

   Recall that ${\cal S}_p$ is the unit sphere in $R^p$.  By the uniform distribution on ${\cal S}_p$ we mean the unique probability distribution on ${\cal S}_p$ that is invariant under the action of ${\cal O}_p$ on ${\cal S}_p$ (see Eaton (1989), Chapter 2). 
  
\medskip 
\noindent{\bf Theorem 2.1:}\ \   Suppose that $\Gamma$ has the Haar distribution on ${\cal O}_p^+$ so $\Gamma_{11}\in (-1,1)$. Let $V_1$  and $V_2$  be independent uniformly distributed vectors in ${\cal S}_{p-1}$.  A version of the conditional distribution of $(\Gamma_{21},\Gamma_{12})$ given $\Gamma_{11}$  is
$${\cal L}(\Gamma_{21},\Gamma'_{12}|\Gamma_{11})={\cal L}\bigl(\sqrt{1-\Gamma_{11}^2}\/ V_1, \sqrt{1-\Gamma_{11}^2}\/ V_2\bigr)\eqno(2.3)$$
where $\Gamma_{11}$ is fixed.

\medskip
   In words, Theorem 2.1 means that given $\Gamma_{11}$, the two functions $(1-\Gamma_{11}^2)^{-1/2}\Gamma_{21}$ and $(1-\Gamma_{11}^2)^{-1/2}\Gamma'_{12}$ of $\Gamma$ are independent and uniformly distributed on ${\cal S}_{p-1}$. A proof of this is given in Eaton and Muirhead (2008).
   
\smallskip
   To specify the conditional distribution $(iii)$ in (2.1), let
$$\eqalign{W_1&=(1-\Gamma_{11}^2)^{-1/2}\Gamma_{21},\cr
W_2&=(1-\Gamma_{11}^2)^{-1/2}\Gamma'_{12}.\cr}\eqno(2.4)$$ 
As stated above, $(W_1,W_2)$  is independent of $\Gamma_{11}$. Further, $W_1, W_2$ are iid uniform on ${\cal S}_{p-1}$.  Since $(\Gamma_{11},\Gamma_{21},\Gamma_{12})$ and $(\Gamma_{11},W_1,W_2)$ are bimeasurable functions of each other, conditioning on $(\Gamma_{11},\Gamma_{21},\Gamma_{12})$ or $(\Gamma_{11},W_1,W_2)$ are equivalent. Next, let $h_i\in{\cal O}_{p-1}$ have first column $W_i$, $i=1, 2$, so 
$$h_i\epsilon_1=W_i,\quad i=1,2$$
 where $\epsilon_1$  is the first standard unit vector in $R^{p-1}$. Then fill out the remaining columns of $h_i$, $i=1,2$   according to Proposition A.2 in Eaton and Muirhead (2008) (with $u=\epsilon_1$ and $\nu=W_i$). 
\medskip
\noindent {\bf Theorem 2.2:}\ \   Suppose  $\Gamma$ is Haar distributed on ${\cal O}_p$.  Let $\Delta$ be Haar distributed on ${\cal O}_{p-2}$.  A version of the conditional distribution of $\Gamma_{22}$ given  $(\Gamma_{11},\Gamma_{21},\Gamma_{12})$ is provided by the distribution 
$${\cal L}\Bigl(h_1\pmatrix{-\Gamma_{11}&0\cr
0&\Delta\cr}h_2'\Bigr)\eqno(2.5)$$
with $h_1$ and $h_2$ fixed and specified above.
\medskip
   The proof of this result can be found in Eaton and Muirhead (2008). One consequence of Theorems 2.1 and 2.2 is that the Haar distribution on ${\cal O}_p$  can be generated by the Haar distribution on ${\cal O}_{p-2}$, the uniform distribution on ${\cal S}_{p-1}$ and the density $f(\cdot |p)$. . Indeed, the following algorithm for generating  $\Gamma$ establishes this claim:
\smallskip 
   
   \item{1.} \ First draw $\Gamma_{11}$ from $f(\cdot |p)$ given in (1.5).
   \smallskip
   \item{2.}\  Next draw iid $U_1$ and $U_2$ which are uniform on ${\cal S}_{p-1}$  and set 
$$\Gamma_{21}=\sqrt{1-\Gamma_{11}^2}\/\/U_1,\quad \Gamma'_{12}=\sqrt{1-\Gamma_{11}^2}\/\/U_2.$$
\item{3.}\  Then construct the $(p-1)\times (p-1)$ orthogonal matrices $h_1$ and $h_2$ (with $W_i=U_i$, $i=1,2$) as above.
\smallskip
\item{4.} Draw $\Delta$ uniformly from ${\cal O}_{p-2}$.
\bigskip
\noindent With
$$A_{22}=\pmatrix{-\Gamma_{11}&0\cr
0&\Delta\cr},$$
the random matrix
$$\Gamma=\pmatrix{\Gamma_{11}&\Gamma_{12}\cr
\Gamma_{21}&h_1A_{22}h_2'\cr}$$
is uniform on ${\cal O}_p$.
 This assertion is a direct consequence of Theorems 2.1 and 2.2.
\bigskip
\goodbreak
\noindent {\bf 3.\ \   The North Pole Distribution}
\medskip
\noindent Here we use the results described in Section 2 to give the distribution of $U_k$ for $k=2, 3$. The notation established in (1.7) and (2.5) is used below. In this notation, some algebra shows that 
$$U_2=\Gamma_{11}^2+(1-\Gamma_{11}^2)W_2'W_1\eqno(3.1)$$
and
$$U_3=\Gamma_{11}^3+2\Gamma_{11}(1-\Gamma_{11}^2)W_2'W_1+(1-\Gamma_{11}^2)W_2'\Gamma_{22}W_1.\eqno(3.2)$$
\smallskip

   We begin with a lemma concerning uniform random vectors on ${\cal S}_{p-1}$.
   \medskip
\noindent {\bf Lemma 3.1:}\ \   Suppose $V_1$  and $V_2$ are independent and uniform on ${\cal S}_{p-1}$. Then the random variable 
$$\xi_2=V_1'V_2\eqno(3.3)
$$ has the density $f(\cdot |p-1)$ defined in (1.5).
\medskip
\noindent{\bf Proof:}\ \   Because ${\cal L}(gV_1)={\cal L}(V_1)$ for all $g\in {\cal O}_{p-1}$,
$${\cal L}(\xi_2)={\cal L}\bigl((gV_1)'V_2\bigr),\qquad g\in {\cal O}_{p-1}.\eqno(3.4)$$
Now condition on $V_1$ and pick $g$ so $gV_1=z_0\in{\cal S}_{p-1}$, where $z_o'=(1,0,...,0)$. Then the conditional distribution of $\xi_2$ given $V_1$ is the same as $z_0'V_2$ which is just the first element of the uniform random vector $V_2$. That this element has density $f(\cdot |p-1)$ is easily established. Since this distribution does not depend on $V_1$, the conclusion of the lemma follows. 
\smallskip
   The proof of Theorem 1.1 is now immediate. From Theorem 2.1, the random variables $W_1$  and  $W_2$ are independent of $\Gamma_{11}$  and are iid uniform on ${\cal S}_{p-1}$. Since  $\xi_{1}\equiv \Gamma_{11}$ has density $f(\cdot |p)$  the representation (1.8) now follows from (3.1) and Lemma 3.1.
\smallskip
   To begin the discussion of ${\cal L}(U_3)$, we first establish the following lemma.
 \medskip  
\noindent{\bf Lemma 3.2:}\ \   Fix $\Gamma_{11}\in (-1,1)$ and let $\Delta$ be uniform on ${\cal O}_{p-2}$. Also, let $u$ and $v$ be fixed vectors in ${\cal S}_{p-1}$  with $u_1$ and $v_1$ denoting the first element of $u$ and $v$ respectively. Then
$${\cal L}\Bigl(u'\pmatrix{-\Gamma_{11}&0\cr
0&\Delta\cr}v\Bigr)={\cal L}\bigl(-u_1v_1\Gamma_{11}+\sqrt{(1-u_1^2)(1-v_1^2)}\/\xi_3\bigr).\eqno(3.5)$$
where $\xi_3$ has density $f(\cdot |p-2)$.
\medskip
\noindent {\bf Proof:}\ \   Since $\Delta$ is uniform on ${\cal O}_{p-2}$, ${\cal L}(g\Delta h)={\cal L}(\Delta)$    for all $g,h\in {\cal O}_{p-2}$. This implies that 
$${\cal L}\Bigl(\pmatrix{-\Gamma_{11}&0\cr
0&\Delta\cr}\Bigr)={\cal L}\Bigl(g^*\pmatrix{-\Gamma_{11}&0\cr
0&\Delta\cr}h^*\Bigr)\eqno(3.6)$$
where 
$$g^*=\pmatrix{1&0\cr
0&g\cr},\quad h^*=\pmatrix{1&0\cr
0&h\cr},\quad g,h\in{\cal O}_{p-2}.$$
Therefore, for $g,h\in{\cal O}_{p-2}$,
$${\cal L}\Bigl(u'\pmatrix{-\Gamma_{11}&0\cr
0&\Delta\cr}v\Bigr)={\cal L}\Bigl((g^*u)'\pmatrix{-\Gamma_{11}&0\cr
0&\Delta\cr}(h^*v)\Bigr).\eqno(3.7)$$
It is clear that we can now pick $g^*$ and $h^*$ so that
$$g^*u=\pmatrix{u_1\cr \sqrt{1-u_1^2}\cr 0\cr\vdots\cr 0\cr} \equiv \tilde u,\quad h^*u=\pmatrix{v_1\cr \sqrt{1-v_1^2}\cr 0\cr\vdots\cr 0\cr} \equiv \tilde v.$$
With these choices for  $g^*$ and $h^*$,
$$\eqalign{{\cal L}\Bigl(u'\pmatrix{-\Gamma_{11}&0\cr
0&\Delta\cr}v\Bigr)&={\cal L}\Bigl({\tilde u}'\pmatrix{-\Gamma_{11}&0\cr
0&\Delta\cr}{\tilde v}'\Bigr)\cr &={\cal L}\bigl(-u_1v_1\Gamma_{11}+\sqrt{(1-u_1^2)(1-v_1^2)}\/\Delta_{11}\bigr)},$$
where $\Delta_{11}$ is the (1,1) element of $\Delta$. Since $\Delta$ is uniform on ${\cal O}_{p-2}$, it follows immediately that $\xi_3\equiv \Delta_{11}$ has density $f(\cdot |p-2)$. This completes the proof of Lemma 3.2.

\smallskip

   Next, the conditional distribution of $U_3$ in (3.2), given $(\Gamma_{11},\Gamma_{21},\Gamma_{12})$,  is described in the following lemma. 
   \medskip

\noindent {\bf Lemma 3.3:}\ \   A version of the conditional distribution of $U_3$  given  $(\Gamma_{11},\Gamma_{21},\Gamma_{12})$ is
$${\cal L}\Bigl(\Gamma_{11}^3+2\Gamma_{11}(1-\Gamma_{11}^2)W_2'W_1+(1-\Gamma_{11}^2)\bigl[-(W_2'W_1)^2\Gamma_{11}+(1-(W_2'W_1)^2)\xi_3\bigr]\Bigr)\eqno(3.8)$$
where $\xi_3$ has density $f(\cdot |p-2)$ and all the remaining quantities in (3.8) are fixed.

\medskip
\noindent {\bf Proof:}\ \ In the notation of Theorem 2.2, recall that the matrix $h_i$ has first column $W_i$, $i=1,2$ and a version of the conditional distribution of $\Gamma_{22}$  given $(\Gamma_{11},\Gamma_{21},\Gamma_{12})$ is
$${\cal L}\Bigl(h_1\pmatrix{-\Gamma_{11}&0\cr 0&\Delta\cr}h_2'\Bigr)\eqno(3.9)$$
with $\Delta$ uniform on ${\cal O}_{p-2}$.
Therefore a version of the conditional distribution of $W_2'\Gamma_{22}W_1$,  given  $(\Gamma_{11},\Gamma_{21},\Gamma_{12})$ , is 
$${\cal L}\Bigl(W_2'h_1\pmatrix{-\Gamma_{11}&0\cr 0&\Delta\cr}h_2'W_1\Bigr)\eqno(3.10)$$
with everything fixed except $\Delta$. Now, we apply Lemma 3.2 with $u=h_1'W_2$ and $v=h_2'W_1$. From the definition of $h_1$ and $h_2$,  it is clear that the first element of both $u$ and $v$ is $u_1=v_1=W_2'W_1$.  From this and (3.5), we conclude that (3.10) is just 
$${\cal L}\Bigl(-(W_2'W_1)^2\Gamma_{11}+\bigl(1-(W_2'W_1)^2\bigr)\xi_3\Bigr).$$
This completes the proof of Lemma 3.3.
\medskip

\noindent {\bf Theorem 3.4:}\ \   Let  $\xi_i$ have density $f(\cdot|p+1-i)$  for  $i=1, 2, 3$ where $\xi_1$, $\xi_2$, and $\xi_3$   are mutually independent. Then 
$${\cal L}(U_3)={\cal L}\Bigl(\xi_1^3+2\xi_1(1-\xi_1^2)\xi_2+(1-\xi_1^2)\bigl[-\xi_1\xi_2^2+(1-\xi_2^2)\xi_3\bigr]\Bigr).\eqno(3.11)$$
  
\medskip
           
\noindent{\bf Proof:}\ \   From Lemma 3.3, a version of the conditional distribution of $U_3$  given $(\Gamma_{11},\Gamma_{21},\Gamma_{12})$  is provided by (3.8). In (3.8), $\xi_3$ is independent of $(\Gamma_{11},\Gamma_{21},\Gamma_{12})$. Setting $\xi_1=\Gamma_{11}$ and $\xi_2=W_2'W_1$,  it follows that $\xi_i$  has density $f(\cdot|p+1-i)$, for  $i=1, 2, 3$. From Theorem 2.1,  $\xi_1$ and $\xi_2$ are independent. By construction, $\xi_3$ is independent of  $(\xi_1,\xi_2)$ so joint independence follows. Unconditioning now yields (3.11) as ${\cal L}(U_3)$. This completes the proof.
\medskip
   The results of Theorem 3.4 provide an easy method to simulate ${\cal L}(U_3)$.  One simply draws $\xi_1$, $\xi_2$, and $\xi_3$ independently as indicated and calculates the algebraic expression in these variables given in (3.11). Note that if $Z_1,\ldots,Z_p$ are iid $N(0, 1)$ variables, then the random variable 
   $$\xi={Z_1\over \sqrt{\sum_{i=1}^p Z_i^2}}$$
has the density $f(\cdot |p)$.
\smallskip
Finally, we note that ${\cal L}({\sqrt{p}}\/U_3)\to N(0,1)$ as $p\to\infty$.  We conjecture that a similar result holds for $U_k$, for any fixed $k$. The method of moments technique used  by Krishnapur (2007) in his proof of Lemma 10 could most likely  be used to establish this, but we have not carried out  the details. Our proof for the cases of $k=1,2,3$ does  not appear to generalize.

\bigskip\bigskip
\centerline{\bf References}
\bigskip

\noindent Diaconis, P. (2003). Patterns in eigenvalues. The 70th Josiah Willard Gibbs Lecture. {\it Bulletin of the American Mathematical Society}, Vol. 40, \#2, 155-178.
\smallskip

\noindent Diaconis, P. and Shahshahani, M. (1994). On the eigenvalues of random matrices. {\it Journal of Applied Probability}, Vol. 31, 49-62.

\smallskip

\noindent Eaton, M. L. (1989). {\it Group Invariance Applications in Statistics}. Volume 1 in the Regional Conference Series in Probability and Statistics published by Institute of Mathematical Statistics, Beachwood, Ohio.

\smallskip
\noindent Eaton, M. L. and Muirhead, R. J. (2008). A decomposition result for the Haar distribution on the orthogonal group. See the arXiv posting http://arxiv.org/abs/0807.2598

\smallskip
\noindent James, A. T. (1954). Normal multivariate analysis and the orthogonal group. {\it Annals of Mathematical Statistics}, Vol. 25, 40-75.

\smallskip
\noindent Krishnapur, M. (2007). From random matrices to random analytic functions. See the arXiv posting http://arxiv.org/abs/0711.1378v1

\smallskip
\noindent Marzetta, T. L., Hassibi, B. and Hochwald, B. M. (2002). Structured unitary space-time autocoding constellations. {\it IEEE Transactions on Information Theory}, Vol. 48, No. 4, 942-950.

\smallskip
\noindent Muirhead, R. J. (1982). {\it Aspects of Multivariate Statistical Theory}. Wiley, New York.

\smallskip
\noindent Rains, E. (1997). High powers of random elements of compact Lie groups. {\it Probability Theory and Related Fields}, Vol. 107, 219-241.

\smallskip
\noindent Sloane, N. (1983). Encrypting by random rotations. Technical Memorandum, Bell Laboratories.

\vfill\eject\bye